\documentclass[11pt]{article}
\usepackage{amssymb,amsfonts,amsmath,latexsym,epsf,tikz,url}

\newtheorem{theorem}{Theorem}[section]
\newtheorem{proposition}[theorem]{Proposition}
\newtheorem{observation}[theorem]{Observation}

\newcommand{\proof}{\noindent{\bf Proof. }}
\newcommand{\qed}{\hfill $\square$\medskip}

\begin{document}

\title{On the dominated chromatic number of certain graphs}

\author{
Saeid Alikhani$^{}$\footnote{Corresponding author}
\and
Mohammad R. Piri  
}

\date{\today}

\maketitle

\begin{center}
Department of Mathematics, Yazd University, 89195-741, Yazd, Iran\\
{\tt alikhani@yazd.ac.ir, piri429@gmail.com}
\end{center}


\begin{abstract}
Let $G$ be a simple graph. The dominated coloring of  $G$ is a proper coloring of $G$ such that each color class is dominated by at least one vertex. The minimum number of colors needed for a  dominated coloring of $G$ is called the dominated chromatic number of $G$, denoted by $\chi_{dom}(G)$.  Stability (bondage number) of dominated chromatic number of $G$ is the minimum number of vertices (edges) of $G$ whose removal changes the dominated chromatic number of $G$. In this paper, we study the dominated chromatic number, dominated  stability and dominated bondage number of  certain graphs.
\end{abstract}


\medskip
\noindent{\bf AMS Subj.\ Class.}: 05C25. 

\section{Introduction and definitions}
In this paper, we are concerned with simple graphs, without directed, multiple, or weighted edges, and without  self-loops. Let $G=(V, E)$ be such a graph and $\lambda \in \mathbb{N}$. A mapping $f: V \longrightarrow \{1, 2,... , \lambda\}$ is called a $\lambda$-proper coloring of $G$ if $f(u) \neq f(v)$, whenever the vertices $u$ and $v$ are adjacent in $G$. A color class of this  coloring is a set consisting of all those vertices assigned the same color. If $f$ is a proper coloring of $G$ with the coloring classes $V_1, V_2, ... ,V_{\lambda}$ such that every vertex in $V_i$ has color $i$, sometimes write simply $f= (V_1, V_2, ... ,V_{\lambda})$. The chromatic number $\chi (G)$ of $G$ is the minimum of colors needed in a proper coloring of a graph. The concept of a graph coloring and chromatic number is very well-studied in graph theory.

A dominator coloring of $G$ is a proper coloring of $G$ such that every vertex of $G$ dominates all vertices of at least one color class (possibly its own class), i.e., every vertex of $G$ is adjacent to all vertices of at least one color class. The dominator chromatic number $\chi_{d}(G)$ of $G$ is the minimum number of color classes in a dominator coloring of $G$. The concept of dominator coloring was introduced and studied by Gera, Horton and Rasmussen \cite{e}.
The total dominator coloring, abbreviated TD-coloring studied well. Let $G$ be a graph with no isolated vertex, the total dominator coloring is a proper coloring of $G$ in which each vertex of the graph is adjacent to every vertex of some (other) color class. The total dominator chromatic number, abbreviated TD-chromatic number, $\chi_{d}^{t}(G)$ of $G$ is the minimum number of color classes in a TD-coloring of $G$. For more information see \cite{nima1, nima2}. 

{\it Dominated coloring} of a graph is a proper coloring in which each color class is a dominated by a vertex. The least number of colors needed for a dominated coloring of $G$ is called the dominated chromatic number  (abbreviated dom-chromatic number) of $G$ and denoted by $\chi_{dom}(G)$ (\cite{Choopani,cc}). We call this coloring a dom-coloring, simplicity.
It is easy to see that for a graph $G$ of order $n$ with $\Delta=n-1$, these parameters 
are equal.

 \begin{proposition}
 	If $G$ is  a connected graph of order $n$ and without isolated vertices with $\Delta (G)=n-1$, then $\chi_{dom}(G)=\chi_d^t(G)=\chi (G)$.
 \end{proposition}

 A set $S$ of vertices in $G$ is a dominating set of $G$ if every vertex of $V(G)\setminus S$ is adjacent to some vertex in $S$. The minimum
 cardinality of a dominating set of $G$ is the domination number of $G$ and denoted by $\gamma (G)$.  
 A domination-critical (domination-super critical, respectively) vertex in a graph
 $G$ is a vertex whose removal decreases (increases, respectively) the domination
 number. Bauer et al. \cite{1}  introduced the concept of domination stability in graphs.
 The domination stability, or just $\gamma$-stability, of a graph $G$ is the minimum number
 of vertices whose removal changes the domination number. Motivated by domination stability, the distinguishing stability and distinguishing bondage number of a graph has introduced in \cite{comm}.

\medskip
 In the next section, we study and compute the dom-chromatic number of certain graphs. In Section 3,   we introduce and study stability and bondage number of dominated chromatic number of graphs.

\section{Dominated chromatic number  of certain graph}

In this section, first we consider graphs obtained by the point-attaching of some another graphs and study their dom-chromatic numbers. 
Also we compute the dom-chromatic number of circulant graphs and cactus graphs.
\subsection{Dominated chromatic number of point-attaching graphs}  

\begin{figure}[h]
	\begin{center}
		\begin{minipage}{5cm}
			\hspace{.4cm}
			\includegraphics[width=3cm,height=3.5cm]{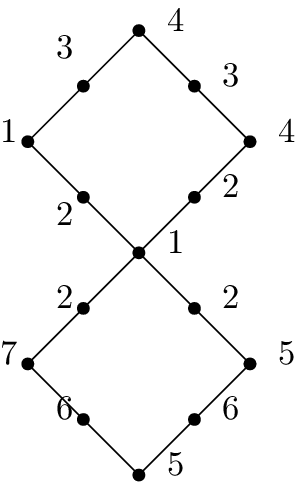}
		\end{minipage}
		\hspace{1.1cm}
		\begin{minipage}{5cm}
			\includegraphics[width=3cm,height=4cm]{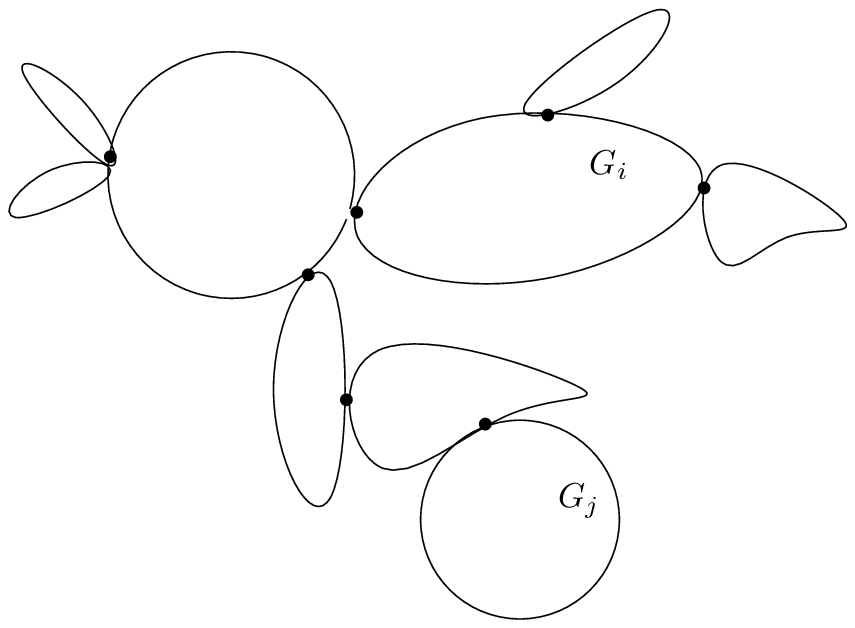}
		\end{minipage}
		\caption{ \label{attaching} \small{ Dom-coloring of  $D_8^2$ and  graph $G$ obtained by point-attaching, respectively.}}
	\end{center}\end{figure}

Let $G$ be a connected graph constructed from pairwise disjoint connected graphs
$G_1,\ldots ,G_k$ as follows. Select a vertex of $G_1$, a vertex of $G_2$, and identify these two vertices. Then continue in this manner inductively.  Note that the graph $G$ constructed in this way has a tree-like structure, the $G_i$'s being its building stones (see Figure \ref{attaching}.  Usually  say that $G$ is obtained by point-attaching from $G_1,\ldots , G_k$ and that $G_i$'s are the primary subgraphs of $G$. A particular case of this construction is the decomposition of a connected graph into blocks (see \cite{Deutsch}). 

In this subsection, we consider some  particular cases  of these graphs  and study dom-chromatic number of these graphs.  The following theorem gives an upper bound for the dom-chromatic number of point-attaching graphs.

\begin{theorem}
	If  $G$ is  point-attaching of graphs $G_1, G_2, ... , G_k$, then 
	\begin{equation*}
	\chi_{dom}(G)\leq \chi_{dom}(G_1)+\chi_{dom}(G_2)+...+\chi_{dom}(G_k).
	\end{equation*}
\end{theorem}
\proof
We color the graph $G_1$ with colors $1, 2, ... , \chi_{dom}(G_1)$, and the graph $G_2$ with colors $\chi_{dom}(G_1)+1, \chi_{dom}(G_1)+2, ... , \chi_{dom}(G_1)+ \chi_{dom}(G_2)$ and do this action for another $G_i$'s. Therefore we have  the result.\qed

\begin{figure}
	\begin{center}
		\includegraphics[width=0.67\textwidth]{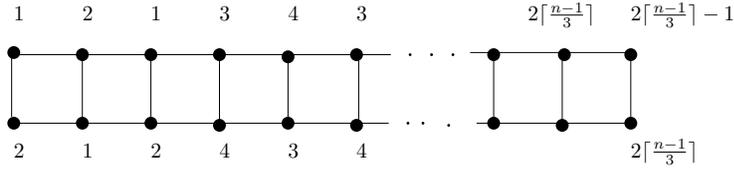}
		\caption{\label{fig1} The dom-coloring of ladder graph $L_n$.}
	\end{center}
\end{figure}

We consider some graphs in the form of point-attaching graphs and compute their dom-chromatic numbers exactly. We need the following theorem: 
\begin{theorem}{\rm(\cite{cc})}\label{path}
For $n\geq 3$,
 \[
 	\chi_{dom}(P_n)=\chi_{dom}(C_n)=\left\{
  	\begin{array}{ll}
  	{\displaystyle
  		\dfrac{n}{2}}&
  		\quad\mbox{if $n\equiv 0 ~(mod\ 4)$, }\\[15pt]
  		{\displaystyle
  			\lfloor \dfrac{n}{2}\rfloor +1} &
  			\quad\mbox{otherwise. }
  				  				\end{array}
  					\right.	
  					\]
\end{theorem}
First  we  consider the ladder graph. Let to recall  the definition of Cartesian product of two graphs. Given any two graphs $G$ and $H$, we define the Cartesian product, denoted $G\square H$, to be the graph with vertex set $V(G)\times V(H)$ and edges between two vertices $(u_1,v_1)$ and $(u_2,v_2)$ if and only if either $u_1=u_2$ and $v_1v_2\in E(H)$ or $u_1u_2 \in E(H)$ and $v_1=v_2$.
The $n-$ladder graph can be defined as $P_2\square P_n$ and denoted by $L_n$. A graph corresponding to the skeleton of an $n$-prism is called a prism graph and is denoted by $Y_n$. A prism graph also called as circular ladder graph, has $2n$ vertices and $3n$ edges. From Figures \ref{fig1} and \ref{prism} which show a dom-coloring of ladder and prism  graphs, we have the following result.

\begin{theorem}
	\begin{enumerate} 
		\item[(i)]  
		For every $n\geq 2$, $\chi_{dom}(L_n)=2\lceil \dfrac{n-1}{3}\rceil$.
		\item[(ii)] 
		For every $n\geq 4$, $\chi_{dom}(Y_n)=\chi_{dom}(L_n)$.
		\end{enumerate} 
				\end{theorem}

\begin{figure}\label{prism}
	\begin{center}
		\includegraphics[width=0.55\textwidth]{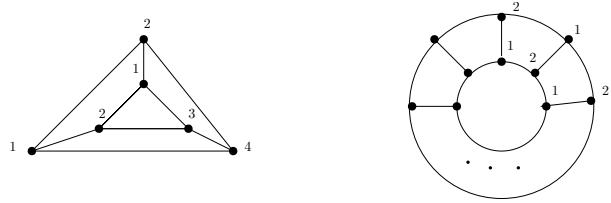}
		\caption{\label{prism}\small{The prism graph.}}
	\end{center}
\end{figure}
\begin{figure}[h]
	\begin{center}
		\begin{minipage}{5cm}
			\hspace{.4cm}
			\includegraphics[width=7cm,height=3.4cm]{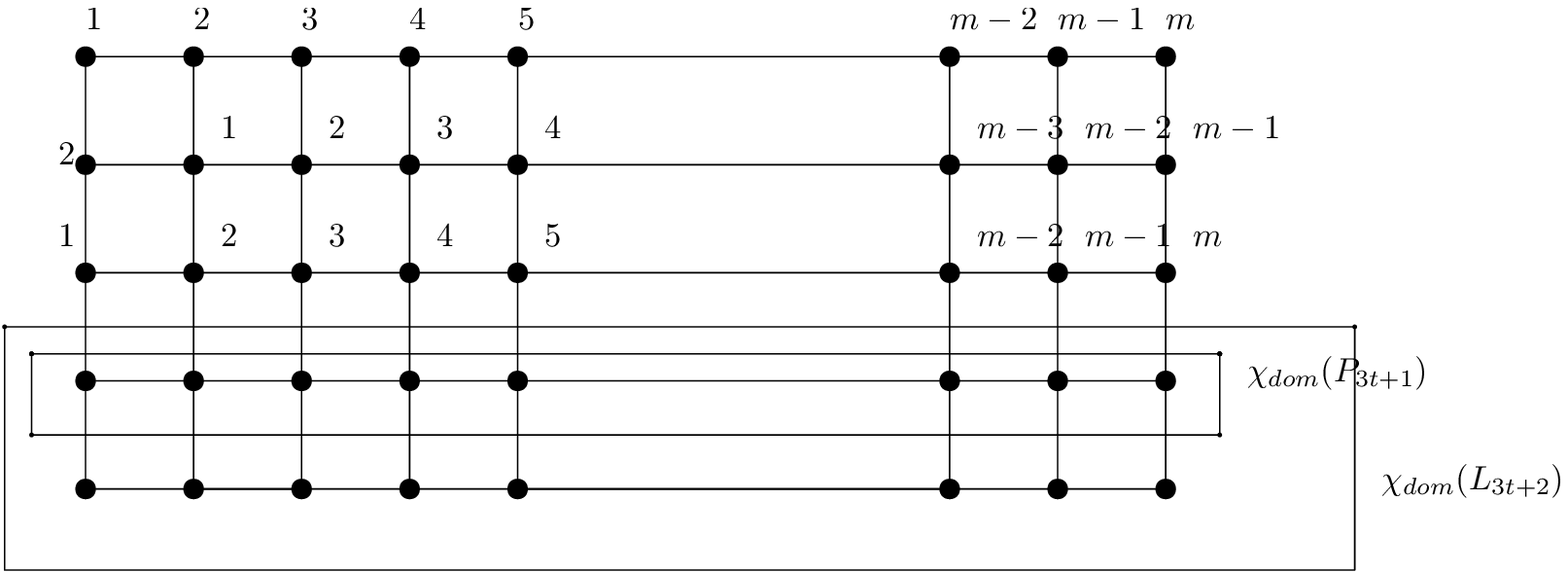}
		\end{minipage}
		\hspace{1.5cm}
		\begin{minipage}{5cm}
			\includegraphics[width=4cm,height=4.5cm]{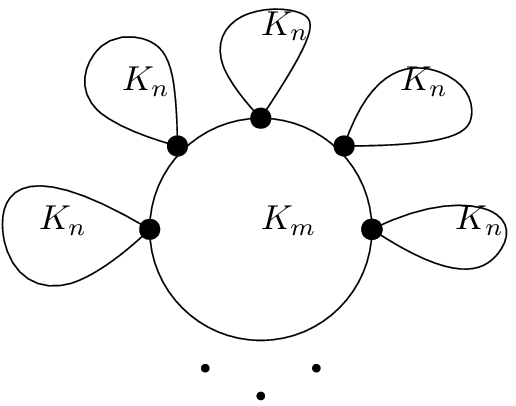}
		\end{minipage}
		\caption{ \label{tour} \small{The grid graph $P_5\square P_m$ and $Q(m,n)$ graph, respectively.}}
	\end{center}\end{figure}

Here, we consider some specific graphs depicted in Figure \ref{tour}.  
Note that the graph $Q(m,n)$ is derived from $K_m$ and $m$ copies of $K_n$ by identifying every vertex of $K_m$ with a vertex of one $K_n$. We compute the dom-chromatic number of the graph $Q(m,n)$. The following theorem gives the dom-chromatic number of grid graphs  
$P_m \square P_n$ and graph $Q(m,n)$. 

\begin{observation}
	\begin{enumerate} 
		\item[(i)]  
		Let $m,n\geq 2$. The dom-chromatic number of gird graphs $P_m \square P_n$ is,
 \[
 	\chi_{dom}(P_n\square P_m)=\left\{
  	\begin{array}{ll}
  	{\displaystyle
  		tm}&
  		\quad\mbox{if $n=3t$, }\\[15pt]
  		{\displaystyle
  			\chi_{dom}(P_m\square P_{3t})+\chi_{dom}(P_m)}&
  			\quad\mbox{if $n=3t+1$,}\\[15pt]
  			{\displaystyle
  				\chi_{dom}(P_m\square P_{3t})+\chi_{dom}(L_m)}&
  				\quad\mbox{if $n=3t+2$.}
  				  					\end{array}
  					\right.	
  					\]
  					\item[(ii)]
  					For $m,n\geq 3$, $\chi_{dom}(Q(m,n))=m(n-1)$.
  					\end{enumerate} 
  					\end{observation}
 \proof
 	\begin{enumerate} 
 		\item[(i)]
 		It follows from the dom-coloring of $P_n\square P_m$ depicted in Figure \ref{tour}.
 		\item[(ii)]
 Let $\lbrace v_1, v_2, ... , v_m\rbrace$ be vertices of $K_m$ and $\lbrace v_j, u_{j_2}, ... , u_{j_n}\rbrace$ be vertices of $j$'th copy of $K_n$ that common in the vertex $v_j$ with $K_m$. We color the vertices of $K_m$ with $\{1, 2,..., m\}$ such that the color $v_i$ be $i$. We color the vertices $u_{j_2}$'s ($j\in\lbrace1, 2, ... , m-1\rbrace$)  with $j+1$ and the vertex $u_{m_2}$ with color $1$ and other vrtices with $\lbrace m+1, m+2, ... \rbrace$. Therefore $\chi_{dom}(Q(m,n))=m+m(n-1)-m=m(n-1)$.\qed
 \end{enumerate}

The friendship graph (or Dutch-Windmill) graph $F_n$ is a graph that can be constructed  by  joining $n$ copies of the cycle graph $C_3$ with a common vertex.  
The generalized friendship graph $D_m^n$ is a collection of $n$ cycles (all af order $m$), meeting at a common vertex. The generalized friendship graph may also be referred to as  flower. We have the following result for the dom-chromatic number of these kind of graphs:
\begin{theorem}
	\begin{enumerate} 
		\item[(i)] 
		For $n\geq 2$, $\chi_{dom}(F_n)=3$.\\
\item[(ii)] 
\[
 	\chi_{dom}(D_m^n)=\left\{
  	\begin{array}{ll}
  	{\displaystyle
  		\chi_{dom}(D_m^{n-1})+\lfloor\dfrac{m}{2}\rfloor}&
  		\quad\mbox{if $n\equiv 1 ~(mod\, 4)$, }\\[15pt]
  		{\displaystyle
  			\chi_{dom}(D_m^{n-1})+\lfloor\dfrac{m}{2}\rfloor -1}&
  			\quad\mbox{otherwise.}
  				  				\end{array}
  					\right.	
  					\]
  					\end{enumerate} 
\end{theorem}
	\proof 
	\begin{enumerate}\item[(i)]
		  Since $F_n$ is the join of $K_1$ and $nK_2$, it is suffices to color 
	$K_1$ with color $1$ and every $K_2$ with colors $2$ and $3$. 
	\item [(ii)]
We color the center vertex with color $1$ and the adjacent vertices  to center  with color $2$ and a vertex adjacent with color $2$  with color $1$. Then for remaining vertices we color as same dominated coloring of  $P_4$ (see Figure \ref{attaching}).\qed
\end{enumerate}

As another example of point-attaching graph,   consider the graph $K_n$ and $n$ copies of connected graph $H$. By definition, the graph $Q(m, H)$ is obtained by identifying each vertex of $K_n$ with a fixed vertex of  $H$.

\begin{theorem}
If the graph $Q(m, H)$ is obtained by identifying each vertex of $K_n$ with a fixed vertex of  $H$, then
\begin{equation*}
n(\chi_{dom}(H)-1)\leq \chi_{dom}(Q(m,H)) \leq n\chi_{dom}(H).
\end{equation*}
\end{theorem}
\proof 
To obtain a dominated coloring for $Q(m,H)$, we have to color the graph $H$ with $\chi_{dom}(H)$ colors and we shall use new colors for another copies of $H$. So we need at most $n\chi_{dom}(H)$ colors. Thus $\chi_{dom}(Q(m,H)) \leq n\chi_{dom}(H)$.\\
Let $V(K_n)=\{v_1, v_2, ... , v_n\}$. Note that $v_i$ is adjacent
 with $v_{i-1}$ and $v_{i+1}$. If we color one of the adjacent vertex of $v_i$ in $H$, with the color of $v_{i-1}$ and do this for another $v_i$'s, then we have dom-coloring and   obviously we need $n(\chi_{dom}(H)-1)$ colors.  So $\chi_{dom}(Q(m,H))\geq n(\chi_{dom}(H)-1)$.\qed

\medskip

Here we consider $r$-gluing of two graphs. 
Let $G_1$ and $G_2$ be two graphs and $r\in \mathbb{N}\cup \{0\}$ with $r\leq min \{\omega(G_1),\omega(G_2)\}$, where $\omega(G)$ shows the clique number of $G$.  Choose a $K_r$ from each $G_i$, $i=1,2$, and form a new graph $G$ from the union of $G_1$ and $G_2$ by identifying the two chosen $K_r$'s in an arbitrary manners. The graph $G$ is called $r$-gluing of $G_1$ and $G_2$ and denoted by $G_1\cup_{K_r}G_2$. If $r=0$ then $G_1\cup_{K_0}G_2$ is just disjoint union. The $G_1\cup_{K_i}G_2$ for $i=1,2$, is called vertex and edge gluing, respectively. 
The following result is about the dominated chromatic number of $r$-gluing of two graphs.

\begin{theorem}\label{glue}
	For any two connected graphs $G_1$ and $G_2$, 
	$$max \{\chi_{dom}(G_1),\chi_{dom}(G_2)\}\leq\chi_{dom}(G_1\cup_{K_r} G_2)\leq\chi_{dom}(G_1)+\chi_{dom}(G_2)-r.$$
\end{theorem}

\proof
Since we need at least $\chi_{dom}(G_1)$ colors to color $G_1$ and $\chi_{dom}(G_2)$ colors to color $G_2$, so we need at least max $\{\chi_{dom}(G_1),\chi_{dom}(G_2)\}$ colors to color $G_1\cup_{K_r} G_2$. So we have $max \{\chi_{dom}(G_1),\chi_{dom}(G_2)\}\leq\chi_{dom}(G_1\cup_{K_r} G_2)$.
\noindent On the other hand, first we give colors $a_1,a_2,\ldots,a_r$ to the vertices of $K_r$. Then we give $a_{r+1},\ldots,a_{\chi_{dom}(G_1)}$ to the other vertices of $G_1$ to have a dom-coloring for $G_1$. Also we give $b_1,b_2,\ldots,b_{\chi_{dom}(G_2)-r}$ to the other vertices of $G_2$ to have a dom-coloring for $G_2$. So every vertex of $G_1\cup_{K_r} G_2$ uses the color class which used before and this is a dom-coloring for $G_1\cup_{K_r} G_2$. So $\chi_{dom}(G_1\cup_{K_r} G_2)\leq\chi_{dom}(G_1)+\chi_{dom}(G_2)-r.$
\qed

\medskip
\noindent{\bf Remark 2.} The bounds in the Theorem \ref{glue} are sharp. For the lower bound it suffices to consider $G_1=K_4$, $G_2=K_5$ and $r=4$. For the upper bound it suffices to consider the complete graph $G_1=K_5$, $G_2=K_6$ and $r=5$.

\subsection{Dom-chromatic number of circulant graphs}  

Let $1\leq a_1<a_2<...<a_m\leq \lfloor \dfrac{n}{2} \rfloor$, where $m, n, a_i$ are integers, $1\leq i \leq m$, and $n\geq 3$. Set $S= \lbrace a_1, a_2, ... , a_m \rbrace$. A graph $G$ with the vertex set $\lbrace 1, 2,..., n\rbrace$ and the edge set 
$\big\{ \lbrace i, j  \rbrace \vert \vert i-j\vert \equiv a_t (mod\, n)$ for some $1\leq t \leq m  \big\}$  
is called circulant graph \cite{f} with respect to set $S$ (or with connection set $S$), and denoted  by $C_n(S)$ or $C_n( a_1, a_2, ... , a_m )$. Notice that $C_n(S)$ is $k$-regular, where $k=2\vert S \vert -1$ if $\dfrac{n}{2}\in S$ and $k=2\vert S \vert $, otherwise. We need the following results to obtain results about dom-chromatic number of circulant graphs.  
\begin{observation}{\rm(\cite{Choopani})}\label{mm}
For any graph $G$, $\chi_{dom}(G)\geq \chi(G)$ and $\chi_{dom}(G)=\chi(G)$ if $diam(G)\leq2$.
\end{observation}
\begin{proposition}{\rm(\cite{cc})}\label{nn}
	\begin{enumerate}
		\item [(i)]
	Let $G$ be a graph with order at least $2$. Then $\chi_{dom}(G)\geq \gamma_t(G)$, where $\gamma_t(G)$ is the total domination number of $G$.
	\item[(ii)] 
	Let $G$ be a graph without isolated vertices. Then $\chi_{dom}(G)\leq \chi(G).\gamma(G)$.
\end{enumerate}
\end{proposition}

From Observation \ref{mm} and Proposition \ref{nn} we have the following proposition. 
\begin{proposition}\label{lll}
If $G$ is  a graph without isolated vertices, then
$$max \lbrace \chi(G), \gamma_t(G)\rbrace \leq \chi_{dom}(G)\leq \chi(G).\gamma(G).$$
\end{proposition}

We need the following theorem:
\begin{theorem}{\rm(\cite{dd})}\label{ddd}
For any $n\geq 4$,
\[
 	\gamma_t(C_n(1,3))=\left\{
  	\begin{array}{ll}
  	{\displaystyle
  		\lceil\dfrac{n}{4}\rceil+1}&
  		\quad\mbox{if $n\equiv 2,4 (mod 8)$, }\\[15pt]
  		{\displaystyle
  			\lceil\dfrac{n}{4}\rceil}&
  			\quad\mbox{otherwise}
  				  				\end{array}
  					\right.	
  					\]
\end{theorem}
\begin{theorem}\label{sss}
If $n=6$, then $\chi_{dom}(C_n(1,3))=3$ and if $n=7$, then $\chi_{dom}(C_n(1,3))=4$ and for $n\geq8$
\[
 	\chi_{dom}(C_n(1,3))=\left\{
  	\begin{array}{ll}
  	{\displaystyle
  		2\lfloor\dfrac{n}{8}\rfloor}&
  		\quad\mbox{if $n\equiv 0 ~ (mod\, 8)$, }\\[15pt]
  	       {\displaystyle
  		2\lfloor\dfrac{n}{8}\rfloor+1}&
  		\quad\mbox{if $n\equiv 1 ~ (mod\, 8)$, }\\[15pt]
  		{\displaystyle
  			2\lfloor\dfrac{n}{8}\rfloor+2}&
  			\quad\mbox{ otherwise.}
  				  				\end{array}
  					\right.	
  					\]
\end{theorem}
\proof
We know $\chi_{dom}(G)\geq \gamma_t(G)$ by Proposition \ref{lll}, and so  
\[
 	\chi_{dom}(C_n(1,3))\geq \gamma_t(C_n(1,3))=\left\{
  	\begin{array}{ll}
  	{\displaystyle
  		\lceil\dfrac{n}{4}\rceil}+1&
  		\quad\mbox{if $n\equiv 2,4 ~(mod\, 8)$, }\\[15pt]
  		{\displaystyle
  			\lceil\dfrac{n}{4}\rceil}&
  			\quad\mbox{otherwise}
  				  				\end{array}
  					\right.	
  					\]
  	Note that $C_6(1,3)$ is isomorphic to complete bipartite graph $K_{3,3}$ and so $\chi_{dom}(C_6(1,3))=\chi_{dom}(K_{3,3})=2$. If $n=7$, then the color classes are: $\lbrace 1,3\rbrace, \lbrace 2,7\rbrace, \lbrace 4,6\rbrace, \lbrace 5\rbrace$. For $n\geq8$ and $1\leq i \leq \lfloor \dfrac{n}{8}\rfloor$ and for $k\in \mathbb{N}$ the color classes are:\\
  	$n=8k$: $A=\lbrace\lbrace 1, 3, n-1, n-3\rbrace \cup \lbrace 2, 4, n, n-2 \rbrace \cup \lbrace 8i-2, 8i, 8i+2, 8i+4\rbrace \cup \lbrace 8i-1, 8i+1, 8i+3, 8i+5 \rbrace \rbrace $\\
  	$n=8k+1$: $A \cup \lbrace n-4 \rbrace$\\
  	$n=8k+2$: $A \cup \lbrace n-4 \rbrace \cup \lbrace n-5 \rbrace $\\
  	$n=8k+3$: $A \cup \lbrace n-4, n-6 \rbrace \cup \lbrace n-5 \rbrace $\\
  	$n=8k+4$: $A \cup \lbrace n-4, n-6 \rbrace \cup \lbrace n-5, n-7 \rbrace $\\
  	$n=8k+5$: $A \cup \lbrace n-4, n-6, n-8 \rbrace \cup \lbrace n-5, n-7 \rbrace $\\
  	$n=8k+6$: $A \cup \lbrace n-4, n-6, n-8 \rbrace \cup \lbrace n-5, n-7, n-9 \rbrace $\\
  	$n=8k+7$: $A \cup \lbrace n-4, n-6, n-8, n-10 \rbrace \cup \lbrace n-5, n-7, n-9 \rbrace $.
  	So we have the result. \qed
  	\begin{theorem}{\rm(\cite{heu})}\label{fdf}
  	Let $C_n(a,b)$ be a circulant graph and $gcd(a,n)=1$. Then the graph $C_n(a,b)$ is isomorphic to graph $C_n(1,c)$ where $c\equiv  a^{-1}b ~ (mod  n)$.
  	\end{theorem}
  	\begin{theorem}
  	Given the circulant graph $C_n(a,b)$, where $n\geq 8$, $gcd(a,n)=1$ and $a^{-1}b\equiv 3 (mod n)$, we have 
  \[
 	\chi_{dom}(C_n(a,b))=\left\{
  	\begin{array}{ll}
  	{\displaystyle
  		2\lfloor\dfrac{n}{8}\rfloor}&
  		\quad\mbox{if $n\equiv 0 (mod 8)$, }\\[15pt]
  	       {\displaystyle
  		2\lfloor\dfrac{n}{8}\rfloor+1}&
  		\quad\mbox{if $n\equiv 1 (mod 8)$, }\\[15pt]
  		{\displaystyle
  			2\lfloor\dfrac{n}{8}\rfloor+2}&
  			\quad\mbox{ otherwise.}
  				  				\end{array}
  					\right.	
  					\]
  	\end{theorem}
  	\proof
  	It follows from  Theorems \ref{sss} and \ref{fdf}.\qed

 \subsection{Dom-chromatic number of cactus graphs} 
 
\begin{figure}[h]
	\begin{center}
		\begin{minipage}{5cm}
			\hspace{.4cm}
			\includegraphics[width=6cm,height=1.4cm]{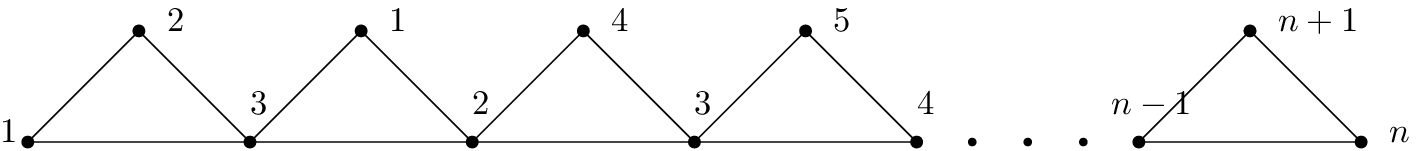}
		\end{minipage}
		\hspace{1.1cm}
		\begin{minipage}{5cm}
			\includegraphics[width=6cm,height=2cm]{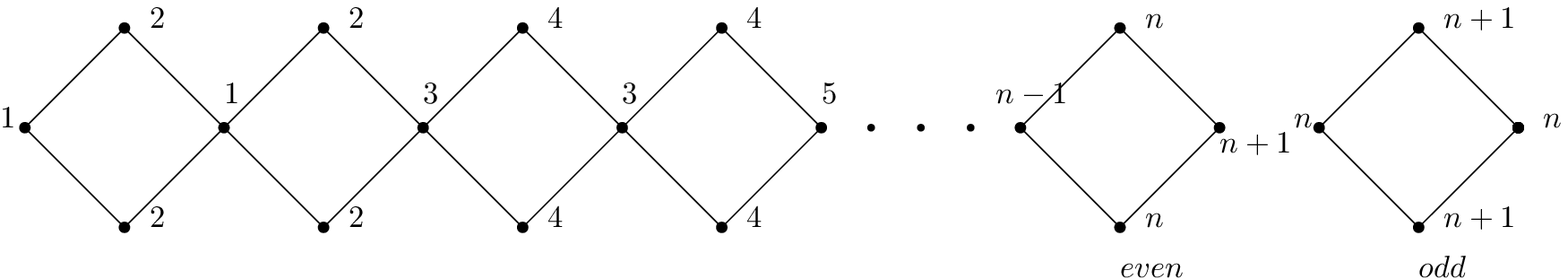}
		\end{minipage}
		\caption{ \label{accd}\small{ The dom-coloring of $T_n$ and $Q_n$, respectively.}}
	\end{center}\end{figure}

   We consider a class of simple linear polymers called cactus chains (\cite{gutman}), in this subsection. A cactus graph is a connected graph in which no edge lies in more than one
cycle. Consequently, each block of a cactus graph is either an edge or a cycle. If all blocks of a cactus $G$ are cycles of the same size $i$, the cactus is $i$-uniform. A triangular cactus is a graph whose
blocks are triangles, i.e., a $3$-uniform cactus. A vertex shared by two or more triangles is called a cut-vertex. If each triangle of a triangular cactus $G$ has at most two cut-vertices, and each cut-vertex is shared by exactly two triangles, we say that G is a chain triangular cactus. By replacing
triangles in this definitions with cycles of length $4$ we obtain cacti whose every block is $C_4$. We call such cacti square cacti. Note that the internal squares may differ in the way they connect to their neighbors. If their cut-vertices are adjacent, we say that such a square is an ortho-square; if the cut-vertices are not adjacent, we call the square a para-square. First we consider a chain triangular. An example of a chain triangular cactus is shown in Figure \ref{accd}. We call the number of triangles in G, the length of the chain. Obviously, all chain triangular cacti of the same length are isomorphic. Hence, we denote the chain triangular cactus of length $n$ by $T_n$   (\cite{gutman}).

\begin{figure}[h]
	\begin{center}
		\begin{minipage}{5cm}
			\hspace{.4cm}
			\includegraphics[width=5.5cm,height=2.3cm]{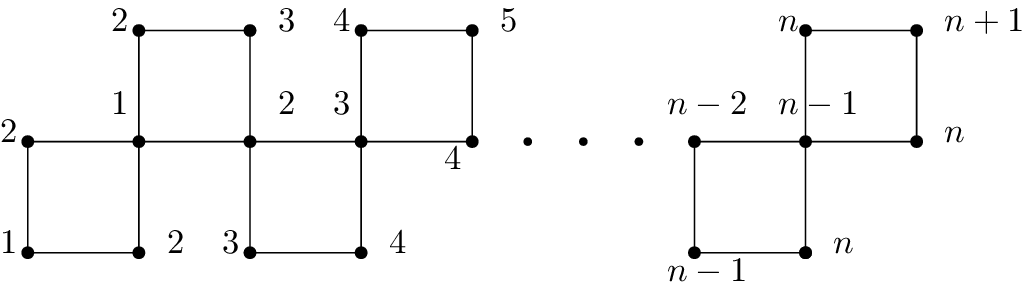}
		\end{minipage}
		\hspace{1.1cm}
		\begin{minipage}{5cm}
			\includegraphics[width=5.4cm,height=2cm]{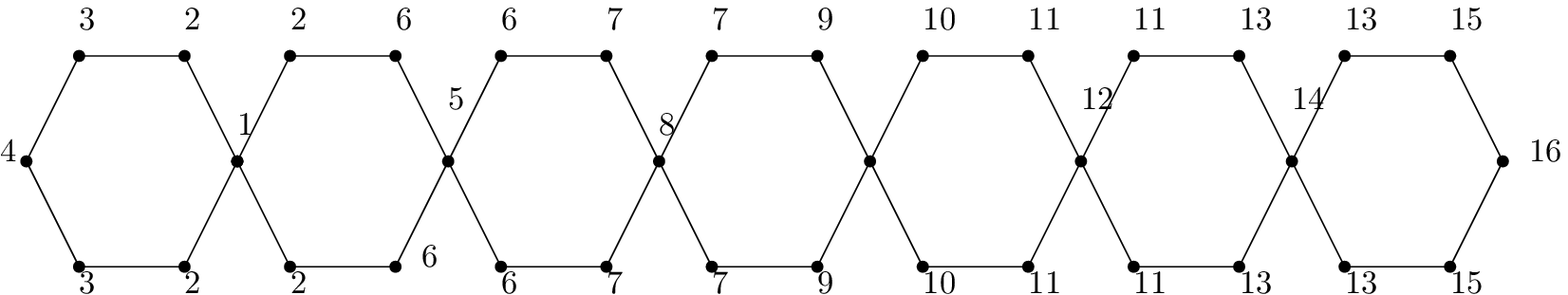}
		\end{minipage}
		\caption{ \label{acdd} \small{The dom-coloring of $O_n$ and $H_n$, respectively.}}
	\end{center}\end{figure}

By replacing triangles in the definitions of triangular cactus $T_n$, with cycles of length $4$ we obtain cacti whose every block is $C_4$. We call such cacti, square cacti. An example of a square cactus chain is shown in Figure \ref{accd}. We see that the internal squares may differ in the way they connect to their neighbors. If their cut-vertices are adjacent, we say that such a square is an ortho-square; if the cut-vertices are not adjacent, we call the square a para-square. We consider a
para-chain of length $n$, which is denoted by $Q_n$ as shown in Figure \ref{accd} and also  another kind of square cactus chain and compute its dominated chromatic number (Figure \ref{acdd}).

	\begin{theorem}
		\begin{enumerate} 
			\item[(i)] 
					For $n\geq 3$, $\chi_{dom}(T_{n})=n+1$.
	\item[(ii)] 	
For $n\geq 2$, $\chi_{dom}(O_n)=\chi_{dom}(Q_n)=n+1$.
\end{enumerate} 
\end{theorem}
\proof
	\begin{enumerate} 
		\item[(i)] We can see in Figure \ref{accd}, 	$\chi_{dom}(T_2)=3$ and for $n\geq 3$,  
		$\chi_{dom}(T_n)=\chi_{dom}(T_{n-1})+1$, and so we have the result.
		
			\item[(ii)] We can see in Figure \ref{acdd}, 	$\chi_{dom}(O_1)=2$ and for $n\geq 2$,  
			$\chi_{dom}(O_n)=\chi_{dom}(O_{n-1})+1$, and so we have the result.\qed
\end{enumerate} 
Replacing triangles in the definitions of triangular cactus, by cycles of length $6$ we obtain cacti whose every block is $C_6$. We call such cacti, hexagonal cacti. An example of a hexagonal cactus chain is shown in Figure \ref{acdd}. We see that the internal hexagonal may differ in the way they connect to their neighbors. If their cut-vertices are adjacent, we say that such a square is an ortho-hexagonal; if the cut-vertices are not adjacent, we call the square a para-hexagonal. We consider a
para-chain of length $n$, which is denoted by $H_n$ as shown in Figure \ref{acdd}.  

\begin{figure}
	\begin{center}
		\includegraphics[width=0.4 \textwidth]{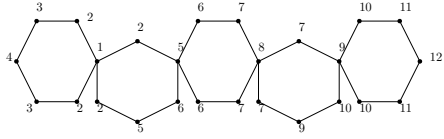}
		\caption{\label{aaccdd1} \small{Meta-chain square cactus graph $M_n$.}}
	\end{center}
\end{figure}

\begin{theorem}
 For every $n\geq 3$, $\chi_{dom}(H_n)=\chi_{dom}(M_n)=n+4$.
\end{theorem}
\proof
By dom-coloring which has shown in Figures \ref{acdd} and \ref{aaccdd1} we have,   
$\chi_{dom}(H_2)=\chi_{dom}(M_2)=6$ and for $n\geq 3$,  
$\chi_{dom}(H_n)=\chi_{dom}(H_{n-1})+2$ and $\chi_{dom}(M_n)=\chi_{dom}(M_{n-1})+2$. 
Therefore we have the result. \qed

\section{Dom-Stability and dom-bandage number of certain graphs}
In this section, we study the stability and bondage number of dominated chromatic number of certain graphs. First we consider stability of certain graphs.
\subsection{Dom-stability of certain graphs}
Stability of dominated chromatic number of a graph $G$, $St_{dom}(G)$, is the minimum number of vertices of $G$ whose removal changes the dom-chromatic  number of $G$. 

\begin{theorem}
For $n\geq 4$
\[
 	St_{dom}(P_n)=\left\{
  	\begin{array}{ll}
  	{\displaystyle
  		2}&
  		\quad\mbox{if $n\equiv 3 ~(mod\, 4)$, }\\[15pt]
  		{\displaystyle
  			1}&
  			\quad\mbox{otherwise}
  				  				\end{array}
  					\right.	
  					\]
\end{theorem}
\proof
Let $V(P_n)=\{v_1,v_2,...,v_n\}$. For $n\geq 4$, we consider the following four cases:
\begin{enumerate} 
\item[(i)]
 If $n=4k$, for some  $k \in \mathbb{N}$, then  in this case by  removing the vertex $v_{4k-1}$, we have  $P_{4k}-v_{4k-1}=P_{4k-2}\cup K_1$. Therefore $\chi_{dom}(P_{4k}-v_{4k-1})=\lfloor \frac{4k-2}{2}\rfloor +1+1 \neq \chi_{dom}(P_{4k})$.

\item[(ii)]  If $n=4k+1$, for some $k \in \mathbb{N}$, then by removing  $v_{4k+1}$, we have $P_{4k+1}-v_{4k+1}=P_{4k}$. Therefore $\chi_{dom}(P_{4+1})\neq \chi_{dom}(P_{4k})$.

\item[(iii)]
 If $n=4k+2$,  for some $k \in \mathbb{N}$, then the proof is similar to proof of Part (i).
 
 \item[(iv)]
  If $n=4k+3$, for some $k \in \mathbb{N}$, then by removing the vertex  $v_{4k+2}$ and  $v_{4k+3}$, we have $P_{4k+3}-\{v_{4k+2}, v_{4k+3} \}=P_{4k+1}$. Therefore $\chi_{dom}(P_{4k+3})\neq \chi_{dom}(P_{4k+1})$.\qed
\end{enumerate} 

\begin{theorem}
For $n\geq 4$, 
\[
 	St_{dom}(C_n)=\left\{
  	\begin{array}{ll}
  	{\displaystyle
  		3}&
  		\quad\mbox{if $n=4k$, }\\[15pt]
  	       {\displaystyle
  		2}&
  		\quad\mbox{if $n=4k+3$, }\\[15pt]
  		{\displaystyle
  			1}&
  			\quad\mbox{ otherwise.}
  				  				\end{array}
  					\right.	
  					\]
\end{theorem}
\proof
Let $V(C_n)=\{v_1,v_2,...,v_n\}$. For $n\geq 4$, we consider the following four cases:
\begin{enumerate} 
	\item[(i)]
 If $n=4k$, for some $k \in \mathbb{N}$, then by removing the vertices  $v_{4k}, v_{4k-1}, v_{4k-2}$, we have  $C_{4k}-\{v_{4k}, v_{4k-1}, v_{4k-2}\}=P_{4k-3}$. Therefore $\chi_{dom}(C_{4k})\neq \chi_{dom}(P_{4k-3})$.
 \item[(ii)]
  If $n=4k+1$, for some $k \in \mathbb{N}$, then by removing the vertex $v_{4k+1}$, we have $C_{4k+1}-\{v_{4k+1}\}=P_{4k}$. Therefore $\chi_{dom}(C_{4k+1})\neq \chi_{dom}(P_{4k})$.
  
  \item[(iii)]
   If $n=4k+2$, for some $k \in \mathbb{N}$, then the proof is similar to the proof of Part (ii).
   \item[(iv)]
    If $n=4k+3$, for some $k \in \mathbb{N}$,  then by removing the vertices $v_{4k+3}, v_{4k+2}$, we have  $C_{4k+3}-\{v_{4k+3},v_{4k+2}\}=P_{4k+1}$. Therefore $\chi_{dom}(C_{4+3})\neq \chi_{dom}(P_{4k+1})$.\qed 
    \end{enumerate}

 The $n$-book graph $(n\geq 2)$ is defined as the Cartesian product $K_{1,n} \square P_2$. We call every $C_4$ in the book graph $B_n$, a page of $B_n$. All pages in $B_n$ have a common side $v_1v_2$. The following observation  gives the dom-stability of $F_n, W_n, D_m^n$ and $B_n$.
\begin{observation}
For $n\geq 2$, $St_{dom}(F_n)=St_{dom}(W_n)=St_{dom}(D_m^n)=St_{dom}(B_n)=1$.
\end{observation}

\begin{figure}
	\begin{center}
		\includegraphics[width=0.24 \textwidth]{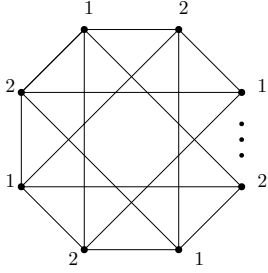}
		\caption{\label{abcd} \small{$n$-regular graph of order $2n$.}}
	\end{center}
\end{figure}

\begin{theorem}\label{egraph}
For every $n\in \mathbb{N}$, there exists a graph $G$ such that $St_{dom}(G)=n$.
\end{theorem}
\proof
Consider the graph $G$ of order $2n$ in Figure \ref{abcd}. As observe that, each vertex with color $1$ is adjacent with color $2$ and each vertex with color $2$ is adjacent to every vertex  with color $1$, and so $\chi_{dom}(G)=2$.
By removing just one vertex of $G$, the coloring  does not change. Suppose that $A$ is the set of vertices whose have color 1. The dom-chromatic number of the induced graph $G-A$ is $n$. The set $A$ has the minimum number of vertices which changes the dom-chromatic number of these kind of graphs (since $K_2$ is always a subgraph of these graphs and we do not need to change the color of the graph by removing each vertex). Therefore $\chi_{dom}(G)=n$.\qed 

\begin{theorem}\label{zzz}
	There exist graphs $G$ and $H$ with the same dom-chromatic number such that $\vert St_{dom}(G)-St_{dom}(H) \vert$ is very large. 
\end{theorem}
\proof 
Suppose that $S_{n_1,n_2}$ is a  double-star with degree sequence $(n_1+1, n_2+1, 1,..., 1)$. It is easy to see that for a bipartite graph  $G=(V_1,V_2,E)$, if there exists a vertex of part $V_1$ which is  adjacent to all vertices of $V_2$ and there exists a vertex in part $V_2$ such that adjacent to all vertices of $V_1$, then $\chi_{dom}(G)=2$. Now 
 let $G=K_{n,n}$ and $H=S_{n_1,n_2}$. So  $\chi_{dom}(K_{n,n})=\chi_{dom}(S_{n_1,n_2})=2$ but $St_{dom}(K_{n,n})=n$ and $St_{dom}(S_{k_1,k_2})=1$. Therefore $\vert St_{dom}(K_{n,n})-St_{dom}(S_{n_1,n_2}) \vert=n-1$.\qed

\begin{proposition}
	If  $G$ is  a graph and $v\in V(G)$, then:
	\begin{equation*}
	St_{dom}(G)\leq St_{dom}(G-v)+1.
	\end{equation*}
\end{proposition}
\proof
If $\chi_{dom}(G)=\chi_{dom}(G-v)$, then $St_{dom}(G)\leq St_{dom}(G-v)+1$ and if $\chi_{dom}(G)\neq \chi_{dom}(G-v)$, then $St_{dom}(G)=1$. So we have the result.\qed

\subsubsection{Dom-bondage of certain graphs}
Bondage number of dominated coloring of graph $G$, $B_{dom}(G)$, is the minimum number of edges of $G$, whose removal changes the dom-chromatic number of $G$. In this subsection, we study the dom-bondage number of specific graphs.
\begin{theorem}
For $n\geq 4$,
\[
 	B_{dom}(P_n)=\left\{
  	\begin{array}{ll}
  	{\displaystyle
  		2}&
  		\quad\mbox{if $n=4k+2$, }\\[15pt]
  		{\displaystyle
  			1}&
  			\quad\mbox{ otherwise.}
  				  				\end{array}
  					\right.	
  					\]
\end{theorem}
\proof
Let $V(P_n)=\{v_1,v_2,...,v_n\}$. For $n\geq 4$, we consider the following four cases:
\begin{enumerate} 
	\item[(i)]
	If  $n=4k$, for some $k\in \mathbb{N}$. In this case, by removing the edge between two vertices $v_{4k-1}$ and $v_{4k}$, we have the result. 
	\item[(ii)]
	 If $n=4k+1$, for some $k\in \mathbb{N}$, by removing the edge between two vertices $v_{4k-1}$ and $v_{4k}$, the new graph is  $P_{4k-1}\cup P_2$ and so  $B_{dom}(P_{4k+1})=1$.
	 \item[(iii)]
	  If $n=4k+2$, for some $k\in \mathbb{N}$, by removing the edges between two vertices $v_{4k-1}$ and $v_{4k+1}$, the new graph is  $P_{4k-1} \cup P_1 \cup P_2$ and so $\chi_{dom}(P_{4k-1} \cup P_1 \cup P_2)\neq \chi_{dom}P_{4k+2}$. 
	  \item[(iv)]
	   If $n=4k+3$, for some $k\in \mathbb{N}$, the proof is similar to the proof of Part (i).\qed
\end{enumerate} 

By removing each edge of the graph $C_n$, we have a path graph of order $n$ and so $\chi_{dom}(C_n)=\chi_{dom}(P_n)$. Therefore we need to remove two edges and so $B_{dom}(C_n)\geq 2$.

\begin{theorem}
For $n\geq 4$,
\[
 	B_{dom}(C_n)=\left\{
  	\begin{array}{ll}
  	{\displaystyle
  		3 }&
  		\quad\mbox{if $n=4k+2$, }\\[15pt]
  		{\displaystyle
  			2} &
  			\quad\mbox{ otherwise.}
  				  				\end{array}
  					\right.	
  					\]
\end{theorem}
\proof
Let $V(C_n)=\{v_1,v_2,...,v_n\}$. For $n\geq 4$, we consider the following four cases:
\begin{enumerate} 
	\item[(i)] If $n=4k$, for some $k\in \mathbb{N}$, then  by  removing the edges between two vertices $v_{4k-1}$ and $v_{v_1}$, the new graph is $P_{4k-1}\cup P_1$. Therefore $B_{dom}(C_{4k})=2$.
	
	\item[(ii)] If $n=4k+1$, for some $k\in \mathbb{N}$, by removing two edges between vertices $v_1$ and $v_2$ and between $v_{4k}$ and $v_{4k+1}$, we have a dom-coloring for new graph with $\lfloor \dfrac{4k+1}{2}\rfloor +2$ colors. So $B_{dom}(C_{4k+1})=2$.
	
	\item[(iii)]
	 If $n=4k+2$, for some $k\in \mathbb{N}$, then by removing  the edges $\{ v_{4k+2}v_1\}$, $\{ v_{4k}v_{4k+1}\}$ and $\{v_{4k-3}v_{4k-2}\}$, the new graph is  $P_{4k-2}\cup 2P_{2}$ and so $B_{dom}(C_{4k+2})=3$.
	 \item[(iv)]
	  If $n=4k+3$,  for some $k\in \mathbb{N}$, the proof is similar to the proof of Part (i).\qed
	  \end{enumerate}

\begin{observation}\label{ss}
For $n\geq 2$, $B_{dom}(F_n)=B_{dom}(B_n)=1$.
\end{observation}

	We end the paper with the following proposition: 

\begin{proposition}
\begin{itemize}
\item[i)]
For any natural number $n\geq2$, there exists a graph $G$ such that $B_{dom}(G)=n$. 

\item[ii)]
There exist graphs $G$ and $H$ with the same dom-chromatic number such that $\mid B_{dom}(G)-B_{dom}(H) \mid$ is very large.
\end{itemize}
\end{proposition}
\proof 
\begin{itemize}
	\item[(i)]
	Since for any $m\geq n$, $B_{dom}(K_{m,n}) = n$, so it suffices to consider $G=K_{m,n}$.
	
	\item[(ii)]
	It suffices  to consider $G=K_{n,n}$ and $H=S_{n_1,n_2}$.\qed
	\end{itemize}


\begin{thebibliography}{99}
 	
 	\bibitem{comm} S. Alikhani, S. Soltani, {\it Stabilizing the distinguishing number of a graph}, Commun. Algebra 46 (2018), no. 12, 5460-–5468. 
	
\bibitem{1} D. Bauer, F. Harary, J. Nieminen, C.  Suffel, {\it Domination alternation sets
in graphs},  Discrete Math. 47  (1983) 153-161.

	
	\bibitem{Choopani}
		F. Choopani, A. Jafarzadeh, D.A. Mojdeh, {\it On dominated coloring of graphs and some Nardhaus-Gaddum-type relations},  Turkish J. Math. 42 (2018) 2148-2156.

\bibitem{Deutsch} E. Deutsch and S. Kla\v{v}zar, {\it  Computing Hosoya polynomials of graphs from primary subgraphs}, MATCH Commun. Math. Comput. Chem. 70 (2013) 627-644.


\bibitem{heu}	C. Heuberger, {\it On planarity and colorability of circulant graphs},  Discrete Math. 268 (2003) 153-169.

\bibitem{Adel} A.P. Kazemi,   {\it Total dominator chromatic number of a graph}, Trans. Combin.  4 (2) (2015)  57-68.

\bibitem{e}
R. Gera, S. Horton, C. Ramussen, {\it Dominator colorings and safe clique partitions}, Congress. Num., 181  (2006) 19-32.



\bibitem{nima1}
N. Ghanbari, S. Alikhani,  {\it More on the total dominator chromatic number of a graph},  J. Inform. Optimiz. Sci.,  40 (2019), no. 1, 157--169.

\bibitem{nima2} 
N. Ghanbari, S. Alikhani,  {\it Total dominator chromatic number of some operations on a graph}, Bull. Comp. Appl. Math., 6 (2018), no. 2, 9-20.

\bibitem{dd} 
N. Jafari Rad,  {\it Domnation in circulant graph}, An. St. Univ. Ovidius Constanta, 17 (2009) 169-176.


\bibitem{cc}
H.B. Merouane, M. Chellali, M. Haddad, H. Kheddouci, {\it  Dominated coloring of graphs}, Graphs Combin. 31 (2015) 713-727.

	\bibitem{f}
	M. Meszka, R. Nedela, A. Rosa,  {\it The chromatic number of $5$-valent circulants},  Discrete Math., 308 (2008) 6269-6284.
	

\bibitem{gutman} A. Sadeghieh, N. Ghanbari, S. Alikhani, {\it Computation of Gutman index of some cactus chains }, Electronic J. Graph Theory  Appl. 6(1) (2018) 138–151. 


\end{thebibliography}
\end{document}